\documentclass[12pt]{amsart}
\input{amssymb.sty}

\newtheorem{theorem}{Theorem}

\newtheorem{lemma}[theorem]{Lemma}

\newtheorem{corollary}[theorem]{Corollary}
\newtheorem{proposition}[theorem]{Proposition}

\begin{document}
\title[Amalgamated products of free groups]
{The Amalgamated product of free groups and residual solvability}
\author[Delaram Kahrobaei]{Delaram Kahrobaei}
\address{ Delaram Kahrobaei, Mathematical Institute, University of St
Andrews, North Haugh, St Andrews, Fife KY16 9SS Scotland, UK}
\email{delaram.kahrobaei@st-andrews.ac.uk}
\urladdr{http://www-groups.mcs.st-and.ac.uk/${\sim}$delaram/}

\begin{abstract} In this paper we study residual solvability of
the amalgamated product of two finitely generated free groups, in
the case of doubles. We find conditions where this kind of
structure is residually solvable, and show that in general this is
not the case. However this kind of structure is always
meta-residually-solvable.
\end{abstract}

\thanks{The research of this author has been supported by CUNY research
foundation at the City College of the City University of New York
and New York Group Theory Cooperative.}

\date{October 2004; submitted}
 \dedicatory{To the memory of Hanna
Neumann and her mathematics} \maketitle
\section{Introduction and Motivation}
The notion of residual properties was first introduced by Philip
Hall in 1954 \cite{PH54}. Let $X$ be a class \cite{HN67} of
groups. $G$ is residually-$X$ if and only if, for every
non-trivial element $g$ in $G$ there is an epimorph of $G$ in $X$
such that the element corresponding to $g$ is not the identity.
There is an important method for constructing new groups from
existing ones called an amalgamated product. These amalgamated
products are special cases of groups acting on trees, embodied in
what is now known as the Bass-Serre theory \cite{S80}.

The question of being residually solvable can be simplified to
whether the group is simple. Peter Neumann in \cite{PN73} asks the
following question: Is it possible that the free product $\{A \ast
B; H = K\}$ with amalgamation, where $A$, $B$ are free groups of
finite ranks, $H, K$ are finitely generated subgroups of $A, B$,
respectively, is a simple group? For the case where the
amalgamated subgroup is not finitely generated, Ruth Camm
\cite{RC53} constructed example of a simple free product $G= \{A
\ast B ; H = K \}$ where $A$, $B$ are free groups of finite rank
and their subgroups $H$, $K$ have infinite rank. This example can
be thought as a simple amalgamated product of two residually
solvable groups where the amalgamation is not residually solvable.
(i.e. where the amalgamated subgroup fails the maximal condition).
In fact she showed that there exist continuously many
non-isomorphic simple amalgamated products of two finitely
generated free groups with non-finitely generated amalgamation.

For the case where amalgamated subgroup is finitely generated
Burger and Mozes \cite{BM97} constructed an infinite family of
torsion free finitely presented simple groups. For every pair $(m,
n)$ of sufficiently large even integers, they constructed a finite
square complex whose universal covering is the product $T_m \times
T_n$ of regular trees of respective degrees $m$ and $n$ and whose
fundamental group ${\Gamma}_{m, n}$ enjoys the following
properties: the group ${\Gamma}_{m, n}$ is simple, finitely
presented and isomorphic to a free amalgamated product $\{F \ast
F; G\}$ where $F$, $G$ are finitely generated free groups.

Note that a partial answer to the question of P. Neumann, is that
under the same condition of the problem, $\{A \ast B; H=K\}$ is
not simple provided either of indices $[A:H], [B:K]$ is infinite;
(see \cite{IS98}, \cite{DK04}, \cite{KS69}).

Clearly, doubles of free groups are never simple. In this paper,
we study residual solvability of amalgamated products of finitely
generated residually solvable groups. Since free groups are
residually solvable, what we prove here is stronger than the
corresponding result for the class of free groups. Indeed we focus
attention on the class of residually solvable groups and show that
the doubles of residually solvable groups are residually solvable
if we impose the condition of solvable separability over the
amalgamating subgroup (see Theorem \ref{RS_solubly_sepr}). In
general we show that doubles of residually solvable groups are not
residually solvable (see Theorem \ref{neg-double})but
meta-residually solvable (see Corollary \ref{meta}).

Fine-Howie-Rosenberger \cite{FHR88} and Culler-Morgan \cite{CM87}
(see also \cite{CGKZ98}) showed that any one-relator group with
torsion which has at least three generators can be decomposed, in
a non-trivial way, as an amalgamated free product. Baumslag-Shalen
\cite{BS} (see also \cite{BS86}) showed that every one-relator
group with at least four generators can be decomposed into a
generalized free product of two groups where the amalgamated
subgroup is proper in one factor and of infinite index in the
other; this is so called Baumslag-Shalen decomposition.

We have recently used these results together with the results of
\cite{DK2004} (see also \cite{DK}) to answer questions regarding
residual solvability of one-relator groups \cite{DK004}. Note that
Baumslag \cite{GB71} proved that positive one-relator groups are
residually solvable (see \cite{KK03} for exposition).

\section{Preliminaries}
In this section, we recall some facts and prove some lemma to be
used later on.
\subsection{Subgroups of amalgamated products}
\label{Hanna} We will use a theorem by Hanna Neumann \cite{HN49}
extensively in this paper. With regard to abstract groups, Hanna
Neumann showed in the 1950s that, in general, subgroups of
amalgamated products are no longer amalgamated products, but
generalized free products, indeed she proved the following: let
$K$ be a subgroup of $G =\{ A \ast B; C \}$, then $K$ is an
HNN-extension of a tree product in which the vertex groups are
conjugates of subgroups of either $A$ or $B$ and the edge groups
are conjugates of subgroups of $C$. The associated subgroups
involved in the HNN-extension are also conjugates of subgroups of
$C$. As a corollary, if $K$ misses the factors $A$ and $B$ (i.e.
$K \cap A = \{1\} = K \cap B$), then $K$ is free; and if $K$
misses the amalgamated subgroup $C$ (i.e. $K \cap C = \{1\}$),
then $K = {\prod_{i \in I}}^* X_i \ast F$, where the $X_i$ are
conjugates of subgroups of $A$ and $B$ and $F$ is free (see
\cite{GB93} for more information).

Let us mention that later a description was given by the
Bass-Serre theory \cite{S80}, with groups acting on graphs to give
geometric intuition: the fundamental group of a graph of groups
generalizes both amalgamated products, HNN-extensions and tree
products.
\subsection{Some Lemmas}
Here we prove some lemmas to be used in proving the main results
of this paper.
\begin{lemma} \label{double_key_lemma} If $A$ is a group, $C$ is a subgroup of
$A$, $\phi$ is an isomorphic mapping of $A$ onto a group $B$, and
$D$ is the amalgamated product of $A$ and $B$ amalgamated $C$ with
$C \phi$, that is
\begin{eqnarray*}
D = \{ A \ast B ; C = C \phi \},
\end{eqnarray*}
then there is a homomorphism, $\psi$, from $D$ onto one of the
factors, and the kernel of $\psi$, $K$, is:
\begin{eqnarray*}
K = gp(a (a \phi)^{-1} | a \in A).
\end{eqnarray*}
Furthermore this map injects to each factor.
\end{lemma}
\begin{proof} Let $\alpha$ be the homomorphism from $A$ onto itself, and $\beta$
be the homomorphism from $B$ onto the inverse of the isomorphic
copy of $A$, i.e. $\beta = {\phi}^{-1}$. These homomorphisms can
be extended to a homomorphism from $D$ onto $A$, (\cite{BN49},
\cite{GB93}). By the way that $\alpha$ and $\beta$ are defined, it
follows that this homomorphism is one-to-one restricted to either
$A$ or $B$.
\end{proof}
\begin{lemma} \label{Not_central_K} Let $A$, $B$, $C$, $D$, $K$ and $\phi$ be as
in Lemma \ref{double_key_lemma}. Then $K$ is not central in $D$,
in other words, $[K, D] \not= \{1\}.$
\end{lemma}
\begin{proof} First, note that the center of $D$ is $\xi D = \xi A \cap \xi B
\cap C.$ By Lemma \ref{double_key_lemma}, $K \cap A = \{ 1 \}$,
and hence $\xi A \cap K = \{ 1 \}$. This implies that $ [K, D]
\not= \{1\}.$
\end{proof}
\begin{lemma} \label{Normal_commute} Let $A$, $B$, $C$, $D$, $K$ and $\phi$ be
as in Lemma \ref{double_key_lemma}. Furthermore, let $C$ be normal
in $A$. Then $K$ commutes with $C$, i.e. $[C, K]=1$.
\end{lemma}
\begin{proof} Let $c, c' \in C$, and $a \in A$. Since $C \lhd A$ and $C \phi
\lhd B$, and $\phi$ is an isomorphism, we can do the following
computation:
\begin{eqnarray*}
a^{-1} c a= c' = c' \phi = (a^{-1} \phi)(c \phi)(a \phi)=(a^{-1}
\phi) c (a \phi).
\end{eqnarray*}
Note that $a (a \phi)^{-1} \in K$, so
\begin{eqnarray*}
c a (a \phi)^{-1} = a (a \phi)^{-1} c;
\end{eqnarray*}
that is $[C, K] = \{1\}$.
\end{proof}
\begin{corollary} Let $A$, $B$, $C$, $D$, $K$ and $\phi$ be as above. Then $K$ is
free.
\end{corollary}
\begin{proof} Since the homomorphism $\psi$ is one-to-one restricted to either $A$ or $B$, then
\begin{eqnarray*}
K \cap A = \{1\} = K \cap A \phi.
\end{eqnarray*}
So, by the theorem of Hanna Neumann mentioned in subsection
\ref{Hanna}, $K$ is free.
\end{proof}

\subsection{The filtration approach to residual solvability}
\label{filtration} In this section we provide some background for
filtration approach which we will use later.

A family $(A_\lambda | \lambda \in \Lambda)$ of normal subgroups
of $A$ is termed a solvable filtration of $A$ if $A/{A_\lambda}$
is solvable for every $\lambda \in \Lambda$ and $\bigcap_{\lambda
\in \Lambda} A_\lambda = \{1\}$. We shall say that $H$ is solvably
separable in $A$ if $\bigcap_{\lambda=1} ^{\infty} H A_\lambda =
H.$ Now let $H \leqslant A$, then $(A_{\lambda}| \lambda \in
\Lambda)$ is called an $H$-filtration of $H$ if $\bigcap_{\lambda
\in \Lambda} H A_{\lambda} = H$. Two equally indexed filtrations
$(A_{\lambda} | \lambda \in \Lambda)$ and $(B_{\lambda} | \lambda
\in \Lambda)$ of $A$ and $B$ respectively are termed $(H,K,
\phi)-$compatible if $(A_\lambda \cap H) \phi = B_\lambda \cap K
\;\; ( \forall \lambda \in \Lambda).$ The following Proposition of
Baumslag \cite{GB63} will help us to prove one of the results: let
$(A_\lambda|\lambda \in \Lambda), (B_\lambda|\lambda \in \Lambda)$
be solvable $(H,K, \phi)-$compatible filtrations of the residually
solvable groups $A$ and $B$ respectively. Suppose $(A_\lambda|
\lambda \in \Lambda)$ is an $H-$filtration of $A$ and $(B_\lambda|
\lambda \in \Lambda)$ is a $K$-filtration of $B$. If, for every
$\lambda \in \Lambda$,
\begin{eqnarray*}
\{ A/{A_\lambda} \ast B/ {B_\lambda};{H A_\lambda}/{A_\lambda}={K
B_\lambda}/{B_\lambda} \},
\end{eqnarray*}
is residually solvable, then so is $G = \{ A \ast B ; H = K \}$.

\section{Doubles of residually solvable groups}
In this section we prove the theorems concerning the doubles of
residually solvable groups. By doubles we mean the amalgamated
product of two groups where the factors are isomorphic and under
the same isomorphism the amalgamated subgroups are identified.

\subsection{Meta-residual-solvability}
Here we prove that in general the amalgamated products of doubles
of residually solvable groups are meta-residually-solvable.
\begin{proposition} \label{easydouble} Let $A$ be a residually solvable group, $C$ be
a subgroup of $A$, and " $\bar{ }$ " be an isomorphic mapping of
$A$ onto $A$. Then the generalized free product of $A$ and
$\bar{A}$ amalgamated $C$ with $\bar{C}$,
\begin{eqnarray*}
G = \{A \ast \bar{A}; C = \bar{C}\}
\end{eqnarray*}
is an extension of a free group by a residually solvable group.
\end{proposition}
\begin{proof} Let $ \phi : G \rightarrow A$; then $K = \ker \phi = gp(a
\bar{a}^{-1}| a \in A)$. $K$ is free by the theorem of Hanna
Neumann mentioned in subsection \ref{Hanna}, since
\begin{eqnarray*}
A \cap K = \{1\} = \bar{A} \cap \ker \phi.
\end{eqnarray*}
Therefore $G$ is an extension of a free group by a residually
solvable group.
\end{proof}
\begin{corollary} \label{meta} $G$ is meta-residually-solvable.
\end{corollary}
\begin{proof}
Since free groups are residually solvable, then by Proposition
\ref{easydouble}, $G$ is residually solvable-by-residually
solvable. That is to say that $G$ is meta-residually-solvable.
\end{proof}
\subsection{Effect of solvably separability on the amalgamated subgroup and residual solvability}
We now prove that if we impose the solvable separability condition
on the amalgamated subgroup of doubles of residual solvable groups
then the resulting group is residually solvable.
\begin{theorem} \label{RS_solubly_sepr} Let " $\bar{  }$ " be an isomorphism from
a group $A$ onto itself, $C$ be a subgroup of $A$, and $G$ be the
amalgamated product of $A$ and $\bar{A}$ amalgamated $C$ with
$\bar{C}$:
\begin{eqnarray*}
G= \{ A \ast \bar{A} ; C = \bar{C}\}.
\end{eqnarray*}
If $A$ is residually solvable then $G$ is also residually
solvable, provided that $C$ is solvably separable in $A$.
\end{theorem}
\begin{proof} Assuming $C$ is solvably separable in $C$ and $A$ is residually
solvable, we want to show that $G$ is residually solvable. That is
we must show that for every non-trivial element $(1 \not=) d \in
G$, there exists a homomorphism, $\phi$, from $G$ onto a solvable
group $S$, $\phi : G \rightarrow S$, such that $d \phi \not= 1.$
We consider two cases:
\newline
Case 1: Let $1 \not= d \in A$. There exists an epimorphism $\phi$
from $G$ onto $A$, so that $d \phi = d$. Since $A$ is residually
solvable, there exists $\lambda \in {\mathbb N}$, such that $d
\not\in \delta_\lambda A$, where $\delta_\lambda A$, is the
$\lambda$-th derived group of $A$. Now put $S=A/{\delta_\lambda
A}$, a solvable group of derived length at most $\lambda$. Note
that the canonical homomorphism, $\theta$, from $A$ onto $S$, maps
$d$ onto a non-trivial element in $S$. Now consider the
composition of these two epimorphism, $\theta \circ \phi$, which
maps $\ G$ onto $S$. The image of $d$ in $S$ is non-trivial:
\begin{eqnarray*}
\theta \circ \phi (d) = d \theta = d \delta_\lambda A {\not=}_S 1.
\end{eqnarray*}
Case 2: Let $1 \not= d \not\in A$ but $d \in G$. Now $d$ can be
expressed as follows:
\begin{eqnarray*}
d = a_1 b_1 a_2 b_2 \cdots a_n b_n \;\; (a_i \in A-C \;\; b_i \in
\bar{A} - \bar{C}).
\end{eqnarray*}
Since the equally indexed filtrations, $\{\delta_\lambda
A\}_{\lambda \in {\mathbb N}}$, and $\{\delta_\lambda \bar{A}
\}_{\lambda \in {\mathbb N}}$ of $A$ and $\bar{A}$ are compatible,
we can form $G_\lambda$:
\begin{eqnarray*}
G_\lambda = \{ A/{\delta_\lambda A} \ast \bar{A}/{\delta_\lambda
\bar{A}} ; {C \delta_\lambda A}/{\delta_\lambda A} = {\bar{C}
\delta_\lambda \bar{A}}/{\delta_\lambda \bar{A}} \}.
\end{eqnarray*}
Note that $G_\lambda$ is residually solvable (by mapping it to one
of the factors and noting that the kernel of the map is free).
Consider the canonical homomorphism $\theta$ from $G$ onto
$G_\lambda$. Since $C$ is solvably separable in $A$, i.e.
\begin{eqnarray*}
\bigcap_{\lambda \in {\mathbb N}} C {\delta}_{\lambda}A & = & C,
\end{eqnarray*}
$\lambda \in {\mathbb N}$ can be so chosen that
\begin{eqnarray*}
a_i \not \in C \delta_\lambda A, \; \bar{a_i} \not\in \bar{C}
\delta_\lambda\bar{A} \; (\text{for } i = 1, \cdots, n).
\end{eqnarray*}
Hence
\begin{eqnarray*}
a_1 \delta_\lambda A \; b_1 \delta_\lambda \bar{A} \cdots a_n
\delta_\lambda A b_n\
 {\delta}_{\lambda} \bar{A} {\not=}_{G_\lambda} 1.
\end{eqnarray*}
This completes the proof of theorem by using Baumslag's
Proposition \cite{GB63}, we recalled in Section \ref{filtration}.
\end{proof}
\subsection{Solvable separability is a sufficient condition for residual
solvability.} Note that the condition of solvable separability of
the amalgamated subgroup in the factors, in the case of doubles,
is necessary. The following theorem shows that the amalgamated
product of doubles is not residually solvable where the factors
are residually solvable groups.
\begin{theorem} \label{neg-double} Let $A$ be a finitely generated residually
solvable group, and $C$ a normal subgroup of $A$, such that $A/C$
is perfect. Let `` $\bar{ }$ `` be an isomorphic mapping of $A$
onto itself. Then
\begin{eqnarray*}
D = \{ A \ast \bar{A} ; C = \bar{C} \}
\end{eqnarray*}
is meta-residually solvable, but {\bf not} residually solvable.
\end{theorem}
\begin{proof} $D$ is meta-residually solvable by Corollary \ref{meta}.
\newline
By Lemma \ref{double_key_lemma}, there exists an epimorphism from
$D$ onto $A$. Let $K$ be the kernel of this epimorphism. Since $C$
is normal in $A$, by using Lemma \ref{Normal_commute}
\begin{eqnarray*}
[C, K]=1. \;\; (*)
\end{eqnarray*}
Now we want to show that $D$ is not residually solvable. We
proceed by contradiction. Suppose $D$ is residually solvable. Let
$d$ be a non-trivial element in $[K,D]$. The existence of such an
element is guaranteed by Lemma \ref{Not_central_K}. Now assume
$\mu$ is a homomorphism of $D$ onto a solvable group $S$, so that
$d \mu \not=_S 1 $. Since $\mu$ is an epimorphism, $C \mu$ is a
normal subgroup of $S$, and by $(*)$,
\begin{eqnarray*}
[C \mu, K \mu ] = 1.
\end{eqnarray*}
If we can show that $S = C \mu$, then $[K \mu, S] = 1$, which
implies that $d \mu =_S 1$, a contradiction. We now need to show
that $S = C \mu$. We have that $D \twoheadrightarrow S$ induces a
homomorphism from $D/C$ onto $S/{C \mu}$. Since $A/C$ and
${\bar{A}}/{\bar{C}}$ each have a perfect subgroup, this induces a
homomorphism from $A/C$ to $1$. So, $S/{C \mu} = 1$ and hence $C
\mu = S$.
\end{proof}
\subsection*{Acknowledgment}
The results of this paper were obtained during my Ph.D. studies at
CUNY Graduate Center and are also contained in my
thesis~\cite{DK}. I express deep gratitude to my Ph.D. supervisor
Gilbert Baumslag whose guidance and support were crucial for the
successful completion of this project. I very much thank Kenneth
Falconer for his sympathetic criticism.
\bibliographystyle{amsplain}
\bibliography{XBib}

\providecommand{\bysame}{\leavevmode\hbox to3em{\hrulefill}\thinspace}
\providecommand{\MR}{\relax\ifhmode\unskip\space\fi MR }
\providecommand{\MRhref}[2]{%
  \href{http://www.ams.org/mathscinet-getitem?mr=#1}{#2}
}
\providecommand{\href}[2]{#2}
\begin{thebibliography}{10}

\bibitem{GB63}
Gilbert Baumslag, \emph{On the finiteness of generalised free products of
  nilpotent groups}, Transaction of American Mathematical Society \textbf{106}
  (1963), 193--209.

\bibitem{GB71}
\bysame, \emph{Positive one-relator groups}, Trans. Amer. Math. Soc.
  \textbf{156} (1971), 165--183.

\bibitem{BS86}
\bysame, \emph{A survey of groups with a single defining relation}, Proceedings
  of {G}roups-{S}t. {A}ndrews 1985, London Math. Soc. Lecture Note Ser.
  Cambridge Univ. Press, Cambridge \textbf{121} (1986), 30--58.

\bibitem{GB93}
\bysame, \emph{Topics in combinatorial group theory}, Birkhauser-Verlag, 1993.

\bibitem{BS}
Gilbert Baumslag and Peter Shalen, \emph{Affine algebraic sets and some
  infinite finitely presented groups}, Essays in group theory, Math. Sci. Res.
  Inst. Publ., 8, Springer, New York (1987), 1--14.

\bibitem{BM97}
Marc Burger and Shahar Mozes, \emph{Finitely presented simple groups and
  products of trees.}, C. R. Acad. Sci. Paris Sér. I Math. \textbf{7} (1997),
  747--752.

\bibitem{RC53}
Ruth Camm, \emph{Simple free products}, J. London Math. Soc. \textbf{28}
  (1953), 66--76.

\bibitem{CGKZ98}
D.J. Collins, R.I. Grigorchuk, P.F. Kurchanov, and H.~Zieschang,
  \emph{Combinatorial group theorey and application to geometry}, Springer,
  1998.

\bibitem{CM87}
Marc Culler and John~W. Morgan, \emph{Group actions on {R}-trees}, Proc. London
  Math. Soc. III \textbf{Ser 55} (1987), 571--604.

\bibitem{FHR88}
Benjamin Fine, James Howie, and Gerhardt Rosenberger, \emph{One-relator
  quotients and free products of cyclics}, Proc. Amer. Math. Soc. \textbf{102}
  (1988), 249--254.

\bibitem{PH54}
Philip Hall, \emph{The splitting properties of relatively free groups}, Proc.
  London Math. Soc. \textbf{4} (1954), 343--356.

\bibitem{IS98}
Sergei~V. Ivanov and Paul~E. Schupp, \emph{A remark on finitely generated
  subgroups of free groups}, Algorithmic problems in groups and
  semigroups(Lincoln NE, 1998),(J.-C. Birget, S. Margolis, J. Meakin, M. Sapir,
  Eds.) (2000), 139--142.

\bibitem{DK2004}
Delaram Kahrobaei, \emph{On residual solvability of gerenalized free products
  of finitely generated nilpotent groups}, preprint (2004), 1--26.

\bibitem{DK004}
\bysame, \emph{On residual solvability of non-positive one-relator groups},
  work in progress (2004).

\bibitem{DK}
\bysame, \emph{Residual solvability of generalized free products}, PhD thesis,
  CUNY Graduate Center (2004), 1--197.

\bibitem{DK04}
\bysame, \emph{A simple proof of a theorem of {Karrass} and {Solitar}}, to
  appear in {A}{M}{S} {P}roceedings {Contemporary} {Mathematics} Geometric
  Method in Group Theory (2004), 1--2.

\bibitem{KK03}
Delaram Kahrobaei and Bilal Khan, \emph{Expository paper of {G}ilbert
  {B}aumslag's theorem: Positive one-relator groups are residually solvable},
  Work in Progress (2004).

\bibitem{KS69}
A.~Karrass and D.~Solitar, \emph{On finitely generated subgroups of a free
  group}, Proc. Amer. Math. Soc. \textbf{22} (1969), 209--213.

\bibitem{BN49}
Bernard~H. Neumann, \emph{An essay on free products of groups with
  amalgamations}, Philos. Trans. Roy. Soc. London. Ser. A. \textbf{246} (1954),
  503--554.

\bibitem{HN49}
Hanna Neumann, \emph{Generalized free producs with amalgamated subgroups. ii},
  Amer. J. Math. \textbf{71} (1949), 491--540.

\bibitem{HN67}
\bysame, \emph{Varieties of groups}, Springer-Verlag New York \textbf{71}
  (1967), 1--192.

\bibitem{PN73}
Peter Neumann, \emph{{S}{Q}-universality of some finitely presented groups}, J.
  Australian Math. Soc. (1973).

\bibitem{S80}
Jean-Pierre Serre, \emph{Trees}, Springer-Verlag, 1980.

\end{thebibliography}

\end{document}